\documentclass{article}
\usepackage{graphicx} 

\usepackage[margin=1in]{geometry}
\usepackage{amsthm, amsmath, amssymb,verbatim,xspace, graphicx}
\usepackage{pdfsync, listings, color,pgf,tikz}
\usepackage[all,cmtip]{xy}
\usepackage{tikz}
\usetikzlibrary{decorations.markings}
\usepackage{amsmath,amssymb,fancyhdr,amsfonts,amsthm,verbatim,graphicx}
\usepackage{latexsym}
\usepackage{exscale}
\usepackage{esvect}
\usepackage{pgf,tikz}
\usepackage{mathrsfs}
\usetikzlibrary{arrows}
\usepackage{cancel}
\usepackage{float}
\usepackage{bm} 
\usepackage{caption}
\usepackage{bigints}
\usepackage{soul}
\usepackage[shortlabels]{enumitem}
\usepackage{ulem}
\usetikzlibrary{trees}

\theoremstyle{definition}
\newtheorem{theorem}{Theorem}[section]

\newtheorem{lemma}[theorem]{Lemma}

\newtheorem{definition}[theorem]{Definition}

\numberwithin{equation}{section}

\def \beq{\begin{equation}}
	\def \eeq{\end{equation}}
\def \lab{\label}
\renewcommand{\rq}[1]{(\ref{#1})}

\renewcommand{\rq}[1]{(\ref{#1})}

\newtheorem{prop}{Proposition}
\newtheorem{thm}{Theorem}
\newtheorem{cor}{Corollary}

\newcommand{\bN}{\Bbb N}

\newcommand{\pa }{\partial }

\newcommand{\Om}{\Omega}
\newcommand{\up}{\Upsilon}

\newcommand{\f}{\varphi }

\newcommand{\Ga}{\Gamma }

\newcommand{\mH}{{\mathcal H}}

\newcommand{\om}{{\omega}}

\def\<{\langle} \def\>{\rangle}

\title{An inverse problem on a metric graph with cycle}
\author{Sergei A. Avdonin and Julian K. Edward}
\date{August 10, 2025}

\begin{document}
	\maketitle
{\bf Abstract} Consider a quantum graph consisting of a ring with two attached edges, and assume Kirchhoff-Neumann conditions hold at the internal vertices. Associated to this graph is a Schr\"{o}dinger type operator $L=-\Delta +q(x)$ with Dirichlet boundary conditions at the two boundary nodes. Let $\{ \omega_n^2, \ \f_n(x)\}$ be the eigenvalues and associated normalized eigenfunctions. Let $v_1$ be a boundary vertex, and $v_2$ the adjacent internal vertex. 
Assume we know the following data: $\{ \omega_n^2,\pa_x \f_n(v_1),\pa_x\f_n(v_2)\}.$ Here $\pa_x\f_n(v_2)$ refers to an outward normal derivative at $v_2$ along one of the edges incident to the other internal vertex.  
From this data we determine the following unknown quantities:  the lengths of edges and the potential functions on each edge. 
\section{Introduction}

\begin{figure}
 \begin{center}
 	\begin{tikzpicture}
 	\begin{scope}[xshift=0cm]
 	\draw (0,0) circle (1cm);
 	\draw[] (-1,0) node[below] {$v_2\quad$} --node[above] {$e_1$}  (-2,0) node[below] {$v_1$}; 
 	\draw[] (330:1cm) node[right] {$v_3$} -- node[below ] {$e_4$} (330:2cm)  node[below] {$v_4$};
 	\draw[] (0,-1) node[above,black] {$e_3$};
 	\draw[] (0,1) node[above] {$e_2$};
 	\fill (-1,0) circle (1.5pt)
 	(-2,0) circle (1.5pt)
 	(-30:1cm) circle (1.5pt)
 	(-30:2cm) circle (1.5pt);;          
 	\end{scope}
 	\end{tikzpicture}
 \end{center}

		\caption{A ring with two attached edges}

\end{figure}
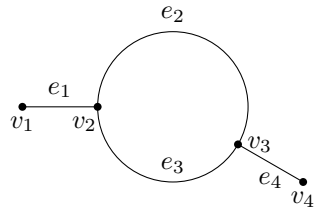

Our inverse problem is considered  on a metric graph consisting of a ring with two edges attached to two different points of the ring, see Figure 1. 
This graph $\Omega =\{ V,E\}$ consists of four vertices,
$V = \{ v_i , i = 1, . . . , 4\}$, and four edges, $E = \{ e_j , j = 1, . . . , 4\}$. We denote the boundary vertices $\{ v_1,v_4\} $ as $\Gamma$.
Denote the length of $e_j$ by $l_j$ . We will often identify edge $e_j$ with the interval $(0,l_j)$, with arclength parametrization $x$. 
We denote by $\phi_j$ the restriction of function $\phi$ to the edge $e_j.$  Assume $q_j$ is a bounded, continuous,
real valued function for all $j$.

Let $\{ (\omega_n^2,\f^n): n\geq 1\}$ be the eigenvalues and normalized eigenfunctions 
of the Laplacian on $\Omega$, with Dirichlet  conditions at the boundary, and Kirchhoff-Neumann (KN) conditions at the interior vertices; thus $\f^n$ solves the following eigenvalue problem:
\beq \lab{sp}
		- \frac{d^2\f^n}{dx^2} + q(x) \f^n = \om_n^2\, \f^n \quad \rm{on} \ \ \{\Omega\setminus V\},
\eeq
\beq \lab{ncf} 
		\sum_{j \in J(v_i)} \pa \f^n_j(v_i)=0, \  v_i \in V \setminus \Ga, \ \f^n|_{\Ga}=0,
  \eeq
\beq \lab{kn}
		\f^n \mbox{ continuous on } \Omega .
\eeq
Here $\f^n_j$ is the restriction of $\f^n$ to edge $e_j$, and $\pa \f^n_j(v_i)$ denotes the derivative of $\f^n$ at $v_i$ pointing away from $v_i$ along $e_j$.

\begin{definition}
The data for the inverse problem, i.e. the items we are assumed to know, are  the spectral data 
$$\{ \omega_n^2,\pa \f^n_1(v_1),\pa \f^n_3(v_2): n\in \bN \}.$$  
\end{definition}
{\bf The inverse problem is then to determine the unknown quantities: $\{ (q_j,l_j):j=1,...,4\}.$}

\begin{thm}\label{thm1}
From the data above, one can recover the set $\{ (l_j,q_j(x)); j=1,...,4\}$.   
\end{thm}

\section{Literature review}

  Inverse problems for differential equations on metric  graphs attracted the attention of the mathematical and physical communities since the 1990s. (The book by Berkolaiko and Kuchment \cite{BerKuch} contains an excellent list of references concerning direct and inverse spectral problem on metric graphs.)
  The first question to be asked when studying inverse problems is how to establish the uniqueness result, i.e. to characterize spectral, or scattering, or dynamical data ensuring a unique solution of the inverse problem. It was shown that inverse boundary spectral and scattering problems for differential equations on graphs with cycles do not, in general, have  a unique solution  \cite{Carl99, KotSmi, GutSmi, KurNow}.
  The results on stable identification are known only for trees, and, almost exclusively for the case of boundary inputs (controls) and observations (inverse problem on trees with internal observations was studied in \cite{AE21}).  
  It was proved that inverse problem on a tree and  is solvable if the actuators and sensors are placed at
  all or all but one of the boundary vertices.

  There are two groups of uniqueness results in this direction: for trees with a priori
  known topology  and lengths of the edges \cite{BrWei, FrYur, Yur05} and for trees with unknown topology \cite{AvKur, AvLeuMik, Bel04, BelVak}.
  The most significant result of \cite{AvKur, AvLeuMik}
  is developing a constructive and robust procedure for the recovery tree's parameters, which became known as the {\bf leaf-peeling  (LP)  method}.
  This method was extended to boundary inverse problems for various types of PDEs on trees
  in a series of our subsequent papers  \cite{AvBell, AvBellNur, AvChLeuMik, AvMikNur, AvNic, AZ21}.
  Our identification procedure is  recursive and  allows
  recalculating efficiently the inverse data from the
  original tree to the smaller trees, ``pruning" leaves step by step up to the rooted
  edge. Because of  its recursive  nature, this procedure  may serve as a basis for developing
  effective numerical algorithms, see, e.g. \cite {AvKhKr, AvKhKr24} as examples of successful numerical realizations. 
  
    The LP method is based on the powerful {\bf boundary control (BC)} method in inverse theory, see \cite{Be87, BeKu, AvBeIv, AvBe96},
  which uses deep connections between controllability and identifiability of distributed parameter systems.
  	The characteristic feature of the BC and LP methods is their {\bf locality}. Specifically,  for inverse problems on graphs,  recovering  the topology and other parameters of a subgraph by the LP method requires only the
  	data related to that subgraph. This property gives the LP method an advantage over other methods and allows us to
  	extend our
  	approach to graphs with cycles. Solving control and inverse problems for differential equations on general graphs requires new developments of the BC and LP methods. 
  In this paper we solve inverse spectral problem for a particular graph with cycle; however, the method we develop can be extended to general graphs, and we plan to do this in our subsequent work.

   Now we describe two results from the literature concerning inverse problems on the graph presented in Figure 1.
   V. Yurko (see \cite{Yur16} Sec. 3.1), considered the spectral problem 
\rq{sp}--\rq{kn} together with three other spectral problems which differs by the boundary conditions: $1)\; \f^n_1(v_1)=0,\, \pa \f^n_4(v_4)=0, \ 2)\; \pa \f^n_1(v_1)=0,\, \f^n_4(v_4)=0, \ 3)\; \pa \f^n_1(v_1)=0,\, \pa \f^n_4(v_4)=0.$ 
Assuming the lengths of all edges are known, he recovered the unknown potentials on all edges from the spectra of these four problems. 

P. Kurasov (see \cite{ Kur23} Sec 23.3), considered the magnetic Schr\"odinger equation. As inverse data he took the Weyl matrix function (i.e. Dirichlet-to Neumann map computed on the graph boundary) known at two values of the magnetic flow. He assumed that spectra of the  Dirichlet operators on the edges forming the cycle have no common points and proved the uniqueness of the inverse problem recovering the potential on the graph.

The statement of our inverse problem and the method of its solution are quite different from these two works. We prove uniqueness for our problem, and our method does not require any conditions on the lengths of the sides. Furthermore, our constructive proof provides an algorithm to recover
 the lengths of the edges and the potential using the minimal possible data (observations).

\section{Solution of the inverse problem}

The statement of our inverse problem is spectral, but its solution will use the dynamical methods. We consider the following initial boundary value problem (IBVP):
\begin{equation} \label{wave1}
u_{tt}-u_{xx} +q(x)u=0 \quad \textrm{ in} \ \  \{\Omega\setminus V\} \times
			(0,T),
\end{equation}
		
\begin{equation} \label{nc}
\sum_{j \in J(v_i)}  \partial u_j(v_i,t)= 0, \ i=2,3,
\end{equation}
\begin{equation} \label{bc}
u_1(v_1,t) = f(t),  u_4(v_4,t)=0,
\end{equation}
\begin{equation} \label{dc}
\begin{cases} 
u_2(v_2,t)-u_1(v_2,t)=0, & \\
u_3(v_2,t)-u_1(v_2,t)=g(t), &   \\    
u_2(v_3,t)=u_3(v_3,t)=u_4(v_3,t), & \end{cases}
\end{equation}	
\begin{equation} \label{ic1}
			u|_{t=0}=u_t|_{t=0}=0 \ \ \textrm{ in} \ \ \Omega.
\end{equation}
We call $f$ and $g$ controls and denote the associated solution by $u=u^{f,g}$. If $g=0$, 
we denote the associated solution by $u^{f,0}$ or simply by $u^{f},$
if $f=0$ -- by $u^{0,g}$ or by $u^{g}.$

We begin with three propositions which are very important for solving our inverse problem. 

\begin{prop}\label{reg}
Let $f,g \in L^2(0,T).$ Then there exists a unique generalized solution of the IBVP \eqref{wave1}--\eqref{ic1}  such that $u^{f,g} \in C(0,T; \mH).$ This solution can be presented in a form of series $u^{f,g}(x,t)=\sum a_n(t)\,\f^n(x)$ where coefficients $a_n(t)$ are UD by the inverse data and controls $f,\,g$:
\beq
 a_n(t)=\int_0^t \left[\pa \f_1^n(v_1)\,f(s)+\pa \f_3^n(v_2) \,g(s)\right]\, \sin \frac{\om_n(t-s)}{\om_n}\,ds
\label{an}\eeq
\end{prop}

The proof of this proposition is based on the Fourier method; the details can be found in \cite{AvNic} Theorems 2,3. 

In Propositions \ref{forward} below we cite our results from \cite{AZ} that are used in the present paper.  This proposition contains several useful formulae regarding solutions of the forward problem for the equation \eqref{wave1} on the simplest graph -- interval $[0,l]$.  

\begin{prop} \label{forward}
Let $w(x,t)$ be the solution to the Goursat problem 
\begin{equation} \label{w}
\begin{cases} w_{tt}-w_{xx}+q(x)w=0, \quad  0<x<t < \infty\\
 w(0,t)=0, w(x,x)=-\frac{1}{2}\int_0^x q(s) \, ds \end{cases}
\end{equation}
in which the $q(x)$ is an extension of the potential function $q(x)$ in \eqref{wave1} from $[0,l]$ to $[0, \infty)$, following the rule $q(2nl \pm x) = q(x)$. Let $k(x,t)$ be the solution to the Goursat problem 
\begin{equation} \label{k}
\begin{cases} k_{tt}-k_{xx}+p(x)k=0,\quad  0<x<t < \infty \\ k( 0,t)=0, \quad k(x,x)=-\frac{1}{2} \int_0^x p(\eta) \, d\eta. \end{cases}
\end{equation}
in which  $p(x)$ is an extension of the potential function $q(x)$ in \eqref{wave1} from $[0,l]$ to $[0, \infty)$, following the rule $p(2nl \pm x)=q(l-x)$ for $0 \le x \le l$.

Let $u^{f,-}$ be the solution to the equation \eqref{wave1} on the interval $[0,l]$ with the initial condition \eqref{ic1}, boundary conditions $u(0,t)=f(t)$ and $u(l,t)=0$;  let $u^{f,+}$ be the solution to the equation \eqref{wave1} on the interval $[0,l]$ with the initial condition \eqref{ic1}, boundary conditions $u(0,t)=0$ and $u(l,t)=f(t)$.  Then: $u^{f,-}$ and $u^{f,+}$ can be expressed in terms of $w$ and $k:$
\begin{multline}  \label{fold} 
u^{f,-}(x,t)   =   f(t-x) + \int_x^t w(x,s)f(t-s) \, ds 
 -f(t-2l+x)-\int_{2l-x}^t w(2l-x,s)f(t-s) \, ds  \\
+f(t-2l-x)+\int_{2l+x}^t w(2l+x,s)f(t-s) \, ds  
-f(t-4l+x)-\int_{4l-x}^t w(4l-x,s)f(t-s) \, ds + \dots  \\
=\sum_{n \geq 0} \left[ f(t-2nl-x) +\int_{2nl+x}^t w(2nl+x,s)\,f(t-s)\,ds \right] \\ + \sum_{n \geq 1} \left[ f(t-2nl+x) +\int_{2nl-x}^t w(2nl-x,s)\,f(t-s)\,ds\right] 
\end{multline}
and \begin{multline}  \label{foldbk} 
u^{f,+}(x,t) = f(t-l+x) + \int_{l-x}^t k(l-x,s) f(t-s) \, ds 
- f(t-l-x) -\int_{l+x}^t k(l+x,s) f(t-s) \, ds +\dots r\\
+f(t-3l+x) + \int_{3l-x}^t k(3l-x,s) f(t-s) \, ds  
- f(t-3l-x) -\int_{3l+x}^t k(3l+x,s) f(t-s) \, ds +\dots \\
=\sum_{n \geq 0} \left[ f(t-(2n+1)l+x) +\int_{(2n+1)l-x}^t k((2n+1)l-x,s)\,f(t-s)\,ds \right] \\ + \sum_{n \geq 0} \left[ f(t-(2n+1)l-x) +\int_{(2n+1)l+x}^t w((2n+1)l+x,s)\,f(t-s)\,ds\right] .
\end{multline}
In \rq{fold},  \rq{foldbk} and everywhere below we assume that our controls $f,\,g$ are extended by $0$ to the negative semiaxis; therefore, for all $t$ the sums above are finite.
\end{prop}
We notice that if $q \in C[0,l]$, then the kernels $w$ and $k$ are continuously differentiable \cite{AvMik}.  

\vskip2mm
Let $u^{f,-}$ be the solution to the equation \eqref{wave1} on the interval $[0,l]$ introduced in Proposition \ref{forward}.
The following well known result is a fundamental tool in our construction:
\begin{prop}\label{basic}
Let $T > l.$
Then the {\it response} operator $R$ defined on $L^2(0,2T)$
by
$$
(Rf)(t)= u^{f,-}_x(0,t), \ t \in (0,2T),
$$
with the domain $\{f \in H^1(0,2T):\; f(0)=0\}$
uniquely determines $l,q.$
\end{prop}
\noindent {\bf Remark 1.}
The same conclusions hold if the edge ${e},$ identified with the interval $(0,l),$ is a part of a graph with the KN conditions holding at $x=0,l$ if we know the operator
$$ u(0,t) \mapsto u_x(0,t), \ t \in (0,2T), \ T>l,$$
defined on $L^2(0,2T).$ This result, based
on the locality of the BC method, was proved in \cite{AvKur},
see also \cite{AZ21}.

\vskip1mm

We now give a constructive proof of Theorem \ref{thm1}. 
In what follows, we denote by UD any data which is uniquely determined, either by hypothesis or by our argument.

\vskip2mm
{\bf Step 1:} Proposition \ref{reg} implies that the operator 
$$ f(t) \mapsto \pa u_1^f(v_1,t)), \ t>0, $$ 
can be computed from our spectral data. According to Proposition \ref{basic} and Remark 1,
knowledge of $R_{11}$ allows us to determine $l_1$ and $q_1$ (independently of the other part of the graph
$\Om \setminus \{e_1\}$).

\vskip2mm
{\bf Step 2:} In this step, we still set $g=0,$ i.e. we have continuity of $u^{f,0}$
at the internal vertex $v_2,$ i.e. $u_j(v_2,t)=:u(v_2,t),\; j=1,2,3.$
With new data from Step 1, we can determine the spectral data 
$$\{ \f^n_1(v_2),\pa\f^n_1(v_2)\}, \mbox{ and hence } \{ \pa\f^n_2(v_2)\} .
$$
Then equation \rq{an} implies that the functions  
$$t\mapsto u^f(v_2,t) \ \ {\rm and} \ \
t\mapsto \pa u^f_j(v_2,t), 
\ j=1,2,3,$$ can be computed from 
$f$ and our inverse data: $$u^f_j(v_2,t)=\sum a_n(t)\,\f^n_j(v_2), \ \  \pa u^f_j(v_2,t)=\sum a_n(t)\, \pa \f^n_j(v_2).$$ 

\vskip2mm
{\bf Step 3:}  
We assume without loss of generality that  $l_2\geq l_3.$ 
 At this step we will compute $l_2,\,l_3,\,q_3$ and a part of $q_2.$ 
First, we wish to conclude that the dynamical Dirichlet to Neumann operator at $v_2$ is UD on $e_j$ for each $j=2,3.$ This requires:
\begin{figure}	
	\begin{center}
		\tikzset{every picture/.style={line width=0.75pt}} 
		
		\begin{tikzpicture}[x=0.75pt,y=0.75pt,yscale=-1,xscale=1]

		\draw    (283.85,78.67) -- (356.85,77.67) ;
		\draw [shift={(356.85,77.67)}, rotate = 358.92] [color={rgb, 255:red, 0; green, 0; blue, 0 }  ][fill={rgb, 255:red, 0; green, 0; blue, 0 }  ][line width=0.75]      (0, 0) circle [x radius= 1.5, y radius= 1.5]   ;
		\draw    (356.85,77.67) -- (394.31,113.57) ;
		\draw [shift={(394.31,113.57)}, rotate = 43.33] [color={rgb, 255:red, 0; green, 0; blue, 0 }  ][fill={rgb, 255:red, 0; green, 0; blue, 0 }  ][line width=0.75]      (0, 0) circle [x radius= 1.5, y radius= 1.5]   ;
		\draw    (356.85,77.67) -- (397.31,47.57) ;
		\draw [shift={(397.31,47.57)}, rotate = 326.58] [color={rgb, 255:red, 0; green, 0; blue, 0 }  ][fill={rgb, 255:red, 0; green, 0; blue, 0 }  ][line width=0.75]      (0, 0) circle [x radius= 1.5, y radius= 1.5]   ;
		\draw [shift={(283.85,78.67)}, rotate = 326.58] [color={rgb, 255:red, 0; green, 0; blue, 0 }  ][fill={rgb, 255:red, 0; green, 0; blue, 0 }  ][line width=0.75]      (0, 0) circle [x radius= 1.5, y radius= 1.5]   ;

		\draw (360,75) node [anchor=north west][inner sep=0.75pt]   [align=left]{$\ \ v_2,x=0$} ;
		\draw (384.04,115.50) node [anchor=north west][inner sep=0.75pt]   [align=left] {$x=l_2$};
		\draw (265.85,83.50) node [anchor=north west][inner sep=0.75pt]   [align=left] {$x=l_1$};	
		\draw (390,30) node [anchor=north west][inner sep=0.75pt]   [align=left] {$ x=l_3$};

		\end{tikzpicture}
		\caption{Star with edges indentified with $(0,l_j)$.}
	\end{center}	
\end{figure}
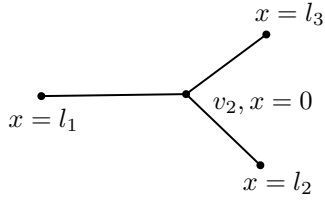
\begin{lemma}\label{onto1}
Let $T>0$. The mapping 
 $
f\mapsto u^f(v_2,t)
$  
is a surjection from $L^2(0,T)$ to $L^2(l_1,T+l_1).$
\end{lemma}
Proof: We will use a formula  proven in \cite{AZ} for general star graphs. For our purposes, the star will consist of edges $e_1,e_2,e_3$ joined at central vertex $v_2,$ see Figure 
2. We identify each edge $e_j$ of the star with the interval $(0,l_j)$, with central vertex $v_2$ identified with $x=0.$
We set $o (t)=u^f(v_2,t)$, and $h_1(t)=u_1^f(v_1,t)=f(t)$, and $h_j(t)=u^f_j(v_3,t)$ for $j=2,3$.  
Assume for the moment that $f\in C^2(0,T)$ with $f(0)=f'(0)=0$, and similar regularity for $h_j,o$, so the solution $u^f$ is classical. Given the functions $o$, $h_2=h_3\in L^2(l_1,T+l_1)$, we will prove the lemma by solving for $h_1=f\in L^2(0,T).$ Using Proposition \ref{forward}, one finds that \rq{nc} implies
\begin{equation} \label{split5}
3 o'(t) - \int_0^t H(s)\, o(t-s) \, ds  =F(t),
\end{equation}
\noindent  where
$$
H(s) = \sum_{j=1}^{3} \partial_x w_j(0, s)
$$
and 
\begin{eqnarray} 
F(t) &= &2 \sum_{j=1}^{3} \sum_{n\geq 0}  h_j'(t-(2n+1)l_j)+2 \sum_{j=1}^{3} \sum_{n\geq 0} k_j( (2n+1)l_j, (2n+1)l_j) h_j(t-(2n+1)l_j) \nonumber \\
&&- 2 \sum_{j=1}^{3}  \sum_{n\geq 0} \int_{(2n+1) l_j}^t \partial_x k_j( (2n+1) l_j,s) h_j(t-s)\, ds -2\sum_{j=1}^{3} \sum_{n\geq 1} o'(t-2nl_j) \nonumber \\
&&-2\sum_{j=1}^{3} \sum_{n\geq 1} w_{j}(2nl_j,2nl_j)o(t-2nl_j)+2\sum_{j =1}^{3} \sum_{n\geq 1} \int_{2nl_j}^t \partial_x w_j (2nl_j, s)o(t-s) \, ds . \label{F5}
\end{eqnarray}
Here $w_j$ and $k_j$ are solutions to \eqref{w} and \eqref{k} on each incident edge of $v_i$.  

Integrating \rq{split5} on $(0,t)$ and using $o (0)=0$, we obtain an integral equation that can be solved iteratively as in \cite{AZ}. For the reader's convenience, we now outline the solution. Our iteration will be on the time intervals $[(2n-1)l_1,(2n+1)l_1)$, with $n\in \bN$.
On the right hand side of the equation below, we will denote by $G(t)$ all terms 
in \rq{F5} that have already been determined in our argument; these
include all terms involving $o ,h_2,h_3$.  
Thus \rq{split5}-\rq{F5}
simplifies, for $t\in [l_1,3l_1)$, to  
$$
G(t)  =2f'(t-l_1)+2k_1( l_1, l_1) f(t-l_1)
-2\int_{l_1}^t \partial_x k_1( l_1,s) f(t-s)\, ds.
$$
Integrating, we get 
$$
\int_{l_1}^tG(s)ds=
2f(t-l_1)+\int_{l_1}^tK(s) f(s-l_1)ds,
$$
with $K\in L^2(l_1,T+l_1).$
Although this equation was derived assuming $f,o,h_j$ were regular, it is standard that  for $h_j,o\in L^2(l_j,l_j+T)$, this Volterra Equation of the Second Kind has a unique solution $f\in L^2(0,T)$, and the resulting $u^f$ will be the weak solution of the wave equation on our star.

Next, we consider \rq{split5}-\rq{F5} for the interval $[3l_1,5l_1)$. Thus, $f(t)$ is known for $t<2l_1$, and again  \rq{split5}-\rq{F5} can be written as 
\beq 
\tilde{G}(t)=
2f(t-l_1)+\int_{3l_1}^tK(s) f(s-l_1)ds,\label{pVESK1}
\eeq
where $\tilde{G}(t)$ is UD.
Arguing as above we thus solve for $f(t)$, $t\in [2l_1,4l_1)$.

Iterating this argument, we solve for $f(t)$ for all $t<T.$
 $\Box$

\ 

The following follows immediately: 
\begin{prop}
Fix $f\in L^2(0,\infty )$.
For $j=2,3$, 
the response operator 
$$u^f_j(v_2,t)\mapsto \pa_x u^f_j(v_2,t),  t>0$$
is UD. 
\end{prop}
\begin{cor}
The following data are UD: $l_2,l_3$, and  $q_3$, and $q_2|_{(0,l_2/2+l_3/2)}.$
\end{cor} 
Proof: All results follow from the remark that follows Proposition \ref{basic} except for the one pertaining to $q_2$. 

We cannot solve for $q_2$ on the entirety of $e_2$ (unless $l_2=l_3$)
 because a wave generated at $v_2$ can pass along $e_3$ to $v_3$, then back to $v_2$ along $e_2$. Upon reaching $v_3$, this wave will contain unknown data from $e_2,e_4$. $\Box$.
 
 To find $q_2$ on the remaining, unknown, portion of $e_2$ will require us to use non-trivial control $g$.

By the corollary, if $l_2=l_3$, then we have uniquely determined $q_2$. So in what follows, we will assume $l_2>l_3.$

\vskip2mm

{\bf Step 4:} We prove that if $g$ is UD on $[0,T]$, then $u^g(v_3,t)$ is UD  on $[0,T]$.

We have proven that $l_3,q_3$ and $\{ \f^n_3(v_2), \pa \f^n_3(v_2)\}$ are UD. By uniqueness of the solution of the Cauchy problem, this proves the data $\{ \f^n_3(v_3), \pa \f^n_3(v_3)\}$ are UD. Therefore, the statement of Step 4  
follows form Proposition \ref{reg}, since $u^{g}(v_3,t)=\sum a_n(t)\,\f^n_3(v_3),$ where 
$$
 a_n(t)=\int_0^t \pa \f_3^n(v_2) \,g(s)\, \sin \frac{\om_n(t-s)}{\om_n}\,ds\,.
$$


\vskip2mm

{\bf Step 5} We use now both controls $f,g$ to find the remaining portion of $q_2$. We will choose $g$ so that the wave generated by $f$ is suppressed along $e_3$, i.e. $u^{f,g}_3(v_2,t)=0$. In effect, this will reduce the wave equation on the graph to a wave equation on a string. The key technical step will then be to ensure that the function $g$ is UD for a given $f$.

We apply \rq{split5}-\rq{F5} in the context of 
 our problem. In particular, we will restrict $u=u^{f,g}$ to edges $e_1,e_2,e_3$, so our domain becomes a three point star. We identify the edge $e_j$ with $(0,l_j)$, with $x=0$ corresponding to the central vertex.
 Thus $u$ will be uniquely determined by $g$ along with $h_j(t):=u(l_j,t)$, $j=1,2,3.$ From the inherited stucture on $\Omega$, we have  $h_1(t)=f(t)$ and $ h_2(t)=h_3(t)$. When referring to the wave on the star, we will drop the superscripts.

\begin{prop}\label{p3}
Let $T>0$, and  $g,h_1,h_2,h_3 \in L^2(0,T)$.
Let $u$ be the solution of the associated wave equation, with zero initial conditions.
Denote $o(t)=u_1(0,t)$. Then

A) the mapping $g \mapsto z(t):=o(t)+g(t)$ is a surjection
$L^2(0,T)\mapsto L^2(0,T).$

B) the mapping $z(t) \mapsto g(t)$ is a bounded mapping from
$L^2(0,T)$ to $L^2(0,T).$
\end{prop}
Proof: 
Fix $z\in L^2(0,T)$ and $h_1,h_2,h_3 \in L^2(0,T)$. 
Then we use \rq{split5}-\rq{F5}  to get
\begin{equation} \label{split3}
3 z'(t) -2g'(t) - \int_0^t \sum_{j=1}^{3} \partial_x w_j(0, s)
 z(t-s) \, ds +\int_0^t \sum_{j=1}^{2} \partial_x w_j(0, s)
g(t-s) \, ds =2F_1(t),
\end{equation}
with
\begin{eqnarray} \label{F3}
F_1(t) &= &  G(t)-\sum_{n=1}^{\lfloor{\frac{t}{2l_j}}\rfloor} (z'(t-2nl_j) )+2\sum_{n=1}^{\lfloor{\frac{t}{2l_j}}\rfloor} g'(t-2nl_j) \nonumber \\
&&-\sum_{j=1}^{3} \sum_{n=1}^{\lfloor{\frac{t}{2l_j}}\rfloor}  w_{j}(2nl_j,2nl_j)z(t-2nl_j)
+\sum_1^2\sum_{n=1}^{\lfloor{\frac{t}{2l_j}}\rfloor}  w_{j}(2nl_j,2nl_j)g(t-2nl_j)\nonumber \\
&&+\sum_{j=1}^{3} \sum_{n=1}^{\lfloor{\frac{t}{2l_j}}\rfloor}  \int_{2nl_j}^t \partial_x w_j (2nl_j, s)z(t-s)-
\sum_{j=1}^{2} \sum_{n=1}^{\lfloor{\frac{t}{2l_j}}\rfloor}  \int_{2nl_j}^t \partial_x w_j (2nl_j, s)g(t-s)
\, ds .
\end{eqnarray}
Here 
 $G(t)\in L^2$ depends only on  $h_1,h_2,h_3$. One can integrate \rq{split3} and then use an iterative Volterra equation argument as in \cite{AZ}  to solve for $g$. 

The proof of B, which is similar, is left to the reader.
 
 $\Box$

We now use the fact that $q_2(x)$ is UD for  $x<l_2/2+l_3/2$ to prove that $q_2(x)$ is UD on $x<3l_2/4+3l_3/4.$ 
\begin{prop}\label{p4}
Let $f\in L^2(0,T)$, and suppose $u^{f,0}$ solves \rq{wave1}-\rq{ic1} on $\Omega$. 
Let $g$ solve
\beq 
u^{0,g}_3(v_2,t)=-u_3^{f,0}(v_2,t),\label{match}
\eeq 
where $u^{0,g}$ solves \rq{wave1}-\rq{ic1} on $\Omega$. 
If $q_2{(x)}$ is UD on the interval $(0, \tilde{l})$, with $\tilde{l}\leq l_2$, then $g(t)$ is UD for $t\in (0, \tilde{l}).$  
\end{prop}
We remark that $g$ solving \rq{match} exists by Proposition \ref{p3}.

\vskip1mm
\noindent{Proof: }

It will be convenient in this proof to set the initial time for $u^{f,g}(x,t)$ at $t=-l_1$, so that the control $g(t)$ will be activated at $t=0.$

Examination of the solution of the Goursat problem (in particular, see Eq.3.5 in \cite{AvMik}) shows that  $q_2(x)$ being UD for $x\in (0, \tilde{l})$, implies that $w_2(x,t)$ is  UD for $x+ t< 2\tilde{l}$. Since $0<x<t$, this implies $w_2(x,t)$ is UD for $x,t<\tilde{l}.$ 

We apply \rq{split3}, \rq{F3}.
In this context, the functions $h_1(t)=0$ is UD for all $t$, whereas
 $h_2(t)=h_3(t)=u^g(v_3,t)$ must be determined by $g$. By unit speed of propagation, we have $h_2(t)=0$ for $t<l_3 $.

Case 1: $\tilde{l}<2l_3.$

Then for $t<\tilde{l}$, the wave generated by $g$ will not have reflected back to $v_1$ along $e_2$ or $e_3$.
 Then, in examining \rq{split3} and \rq{F3}, and labelling various UD functions by $G$, we have 
$$
G(t)+g'(t)-\int_0^t \sum_j\pa_xw_j(0,s)g(t-s)ds
$$
\beq =
-2\sum_{j\neq 2}\sum_{n\geq 1}
\left ( g'(t-2nl_j)+w_j(2nl_j,2nl_j)g(t-2nl_j)-\int_{2nl_j}^t\pa_xw_j(2nl_j,s)g(t-s)ds
\right ).\label{ve}
\eeq 
Here we use the rule that $g(t)=0$ for $t<0$, so the sum on the right is finite. Observe that all terms on except $g$ are UD, so 
 $g(t)$ is UD for such $t$. Thus in this case,   the proof is complete.

Case 2: $\tilde{l}\geq 2l_3.$

In this case, we first suppose $t<2l_3$.
Then arguing as in Case 1, we can solve for $g(t)$  on $(0,2l_3)$ and it will be UD. Thus, by Step 4,  $h_2(t)=h_3(t)$ are UD on the same time interval. We remark in passing that $h_3(t)=0$ for $t<l_3$ by unit speed of propagation.
  Then, in examining \rq{split3} and \rq{F3}, and labelling various UD functions by $G$, we have that \rq{ve} again holds, now for $(0,\max (3l_3,\tilde{l})$. Thus all terms  except $g(t)$ are UD, so 
 $g(t)$ is UD for such $t$. If $\tilde{l}<3l_3$, then the proof is complete. If we  assume otherwise, then we can iterate this argument until 
 $\tilde{l}<nl_3$ for some $n$, thus completing the proof. 
 $\Box$

\ 

\begin{lemma}\label{L1}
Let $f\in L^2(0,T)$, and let $g\in L^2(l_1,T)$ satisfy \rq{match}. Then the mapping 
$$
f\mapsto u_1^{f,g}(v_2,t)
$$
is  a surjection onto $L^2(l_1,T)$.
\end{lemma}
Proof: 
Define the operators 
$$
Af(t)=u_3^{f,0}(v_2,t),\ Bg(t)=u_3^{0,g}(v_2,t),\ 
Cg(t)=u_1^{0,g}(v_2,t).
$$
Then by \rq{dc} we have $B=C+I$, with $I$ the identity operator, and by \rq{match} we have $Af+Bg=0.$ By Proposition \ref{p3} part B, $B$ is invertible, hence $g=-B^{-1}Af$. Then
\begin{eqnarray*}
 u_1^{f,g}(v_2,t)& = & Af(t)+Cg(t)\\
 &= & Af(t)-CB^{-1}Af(t)\\
 & = & B^{-1}Af(t).
\end{eqnarray*}
By Lemma \ref{onto1}, $A$ is onto, and obviously the same holds for $B^{-1}$, the lemma follows. $\Box$

We now apply Proposition \ref{p4}, Lemma \ref{L1} to complete Step 5. 
We can set $\tilde{l}= l_2/2+l_3/2.$ 
We will use $f$ to start a wave at $v_0$ starting at $t=0$. Then at $t=l_1$, we use $g$, with $g(t)=0$ for $t<l_1$,  so that $u_3^{f,g}(v_2,t)$ will be zero 
for  $t<l_1+\tilde{l}$. Then evidently, $u^{f,g}(v_3,t)=0$ for $t<l_1+\tilde{l}+l_3.$

Now let $T\in (l_1,l_1+{2}l_2)$, and $f\in L^2(0,T).$ Then let $g(t)$ as in Proposition \ref{p4} but with initial time set at $t=l_1$, so $g$ is first 
defined 
for $l_1<t<l_1+\tilde{l}$ . We then extend $g(t)=0$
for $t>l_1+\tilde{l}$ and for $t<l_1$. 
Now
\begin{eqnarray*}
u^{f,g}_3(v_2,t)& =& u^{f,0}_3(v_2,t)+u^{0,g}_3(v_2,t)\\
& = & 0,\ t<l_1+\tilde{l}.
\end{eqnarray*}
Because $u^{f,g}_3(v_2,t)=0$ for $t<l_1+\tilde{l},$ it follows that
a wave from $v_2$ will not travel along $e_3$ to reach $v_3$ in time $t<l_1+\tilde{l}+l_3$.  Thus, any wave travelling first along $e_3$, and then $e_2$,  will only reach  $v_2$  by time
$$
t=l_1+\tilde{l}+l_3+l_2=l_1+3l_2/2+3l_3/2=:l_1+\hat{l}.
$$ 
Thus, the function
$$
t\mapsto (u_2^{f,g}(v_2,t),\pa u_2^{f,g}(v_2,t)), t<l_1+\hat{l}
$$
is UD. By Lemma \ref{onto1}, it follows that 
for the edge $e_2$, which we identify with $(0,l_2)$ with $v_2$ corresponding to $x=0$, {the mapping 
$$
u_2^{f,g}(v_2,t)\mapsto \pa u_2^{f,g}(v_2,t)
$$
is UD for $t< \hat{l}.$
By Lemma \ref{basic}, we conclude
 $q_2(x)$ is UD for $x<\hat{l}/2$.
 If $3l_3>l_2$, then we have  $\hat{l}>l_2$, and hence completely solved for $q_2$.
 Otherwise, it suffices to set $\tilde{l}=3l_2/4+3l_3/4$ and repeat the argument above.

{\bf Step 6:} we solve IP on $e_4.$ In this step, we can set $g=0.$

Having solved for $q_2$ and $l_2,$ we can use uniqueness of  the solution of the Cauchy problem to find $\{ \pa \f^n_2(v_3)\}$, and hence $\{ \pa \f^n_4(v_3)\}$. Hence, the mapping 
$$f\mapsto \pa u_4^{f,0} (v_3,t)$$
is  UD. Furthermore, again by Cauchy uniqueness and the arguments of Lemma \ref{onto1}, the mapping 
$$f\mapsto  u_4^{f,0} (v_3,t)$$
is both UD and a surjection. This shows that the response operator
{the mapping 
$$
u_4^{f,0}(v_3,t)\mapsto \pa u_4^{f,0}(v_3,t)
$$
is UD for $t< \hat{l}.$}
 It follows by Proposition \ref{basic} that $l_4,q_4$ are UD.

 \vskip3mm
 \noindent {\bf  Acknowledgments}\\
The research of Sergei Avdonin was  supported  in part by the National Science Foundation, grant DMS 2308377.

\end{document}